\documentclass[reqno,10pt]{amsart}
\usepackage{amsmath,amsthm,amssymb,amscd}
\advance \topmargin by-\headheight
\advance \topmargin by-\headsep
\evensidemargin \oddsidemargin
\theoremstyle{plain}

\newtheorem{thm*}{Theorem}
\newtheorem{lem*}{Lemma}

\theoremstyle{definition}

\numberwithin{equation}{section}

\newcommand{\CC}{\mathbb{C}}

\begin{document}
\title{\bf{Three Lectures on the Riemann Zeta-Function}}
\author{{\bf S.M. Gonek} \\
Department of Mathematics \\ 
University of Rochester \\
Rochester, N. Y. 14627 \\
U.S.A.}
\thanks{The work of the author was supported in part by a grant from 
the National Science Foundation.}

\date{}

\maketitle

\section*{Introduction}

These lectures were delivered at the ``International 
Conference on Subjects Related to the Clay Problems'' held at Chonbuk National 
University, Chonju, Korea in July, 2002. My aim
was to give mathematicians and graduate students unfamiliar 
with analytic number theory an introduction to the theory 
of the Riemann zeta--function focusing, in particular, on the 
distribution of its zeros. Professor Y. Yildirin of the 
University of Ankara, who also delivered a set of lectures at the 
conference, concentrated on the distribution of prime numbers.

A few general remarks about the lectures are in 
order before I summarize their contents. First, since I could only cover 
a small part of the subject in the time alloted, my choices about what to include and 
exclude were necessarily personal. Second, I have glossed over a number of 
technical details in order to keep the focus on the main ideas. 
Finally, there is almost nothing new in the lectures. The 
exception is the description of a new random matrix model due to 
C. Hughes, J. Keating, and the author at the end of the third lecture.
I should also add that this manuscript is a very close record of the lectures 
I delivered and this, I think, accounts for the somewhat breezy style.

In the first lecture I presented the basic background 
material on the zeta--function, sketched a proof of 
the Prime Number Theorem, explained how the Riemann 
Hypothesis (RH) comes into the picture, and briefly summarized the evidence 
for it.

In the second lecture I wanted to explain how one
studies the distribution of the zeros and chose
mean--value estimates as a unifying theme. 
I described what mean--value estimates are, gave several examples, and 
explained in a general way their connection with the zeros.
I then sketched the ideas behind two applications -- the most primitive 
zero--density estimate (due to H. Bohr and E. Landau) and the 
proof of N. Levinson's famous result that at least one--third of the 
zeros of the zeta--function lie on the critical line. 
Both results were cited in Lecture I as evidence for the 
Riemann Hypothesis. I had also intended to present the conditional result
of J. B. Conrey, A. Ghosh, and the author that more than seventy percent of 
the zeros are simple, but there was not enough time. However, I have included 
that application here.

The third lecture began with the observation that
the Riemann Hypothesis does not answer all our
questions about the primes; one also needs  
detailed information about the vertical distribution of the 
zeros on the critical line. I then presented H. Montgomery's 
pioneering work on the pair correlation of the zeros. In the remainder 
of the lecture I stated
the GUE hypothesis  and described the most recent work on modeling the zeta--function by 
characteristic polynomials of random matrices from the Circular 
Unitary Ensemble (CUE). 

For those wishing to study the zeta--function in more depth, the most 
important books are by H. Davenport [D], H. M. Edwards [E], A. E. Ingham 
[I2], A. Ivic [Iv], and E. C. Titchmarsh [T1], [T2 ]. For a background in random 
matrix theory the reader should consult M. L. Mehta [M] and P. Deift 
[Df].

I take this opportunity to thank the organizers and the many other fine 
Korean mathematicians I got to meet for the first time at the conference. Thanks also to 
the mathematicians and students who so warmly hosted us visiting 
mathematicians and made the conference such an enjoyable and memorable one.
\newpage

\centerline{\bf\LARGE{Lecture I}} 
\medskip  
\centerline{\bf{\Large{The Zeta--Function, Prime Numbers, 
and the Zeros}}}
\medskip

Although most mathematicians are aware that the prime numbers, the Riemann 
zeta--function, and the zeros of the zeta--function are intimately 
connected, very few know why. In this first lecture I will outline the basic 
properties of the zeta--function, sketch a proof of the prime number theorem, 
and show how the location of the zeros of the zeta--function directly influences 
the distribution of the primes. I will then explain why the Riemann Hypothesis (RH) is 
important and the evidence for it.

\noindent \section*{1 \quad The Riemann zeta-function}

The Riemann zeta--function is defined by the Dirichlet series
$$ 
\zeta (s) = \sum^\infty_{n =1} n^{-s} ,
$$
which can also be written
$$ 
\zeta (s) = \prod_p (1 + p^{-s} + p^{-2s} + \cdots) = \prod_p (1 -
 p^{-s})^{-1}\;,
$$
where $s=\sigma+it$ is a complex variable. We immediately see that the zeta--function 
is built out of the prime numbers.
Observe that the series and product both converge absolutely in 
the half--plane $\sigma > 1$. Their equality in this region
may be regarded as an analytic equivalent of the Fundamental Theorem of 
Arithmetic. For the Fundamental Theorem assures us that each term $n^{-s}$ in 
the series occurs once, and only once, among the terms resulting from 
multiplying out the Euler product. Conversely, if we know the 
equality of the sum and product, the Fundamental Theorem follows.
From the equality of the sum and product we can also
deduce the well kown fact that there are an infinite number of primes. 
For if there were not, the product would remain bounded as $\sigma \to 1^{+}$, whereas 
we know that the sum tends to infinity. 

Since no factor in the Euler product equals zero when $\sigma > 1$, we 
deduce that $\zeta (s) \ne 0 $ when $\sigma > 1$ \,. Also, since the series converge 
absolutely when $\sigma > 1$, 
it converges uniformly in compact subsets there. It follows that 
$\zeta (s)$ is  analytic in the half--plane $\sigma > 1$.

The most fundamental properties of  the zeta--function are:
\begin{enumerate}
 
{\item {\bf{ Analytic continuation:}}}
 $\zeta (s)$ has an analytic continuation to $\mathbb{C}$ except
 for a simple pole at $s =1$.
 
 \item {\bf{ Functional equation:}}  The zeta--function satisfies the 
 functional equation  
 $$
 \pi^{-s/2} \Gamma(s/2)\zeta (s) = 
 \pi^{-(1-s)/2} \Gamma((1-s)/2)\zeta (1-s)\,. 
 $$
 
 \item  {\bf{ Trivial zeros:}} The only zeros of $\zeta (s)$ in 
 $\sigma <0$ are simple ones at $s = -2, -4, -6, \dots$\,.
 
 \item {\bf{Nontrivial zeros: }} $\zeta (s)$ has infinitely many
 zeros $\rho = \beta + i \gamma$ in the ``critical strip'' $0 \leq 
 \sigma \leq 1$.  These lie symmetrically about the ``critical line'' 
 $\sigma = 1/2$, and about the real axis. 
 
\item {\bf{Density of zeros in the critical strip: }}
If $N(T)$ denotes the number of zeros $\rho = \beta + i \gamma$ in the 
critical strip with ordinates $0< \gamma \leq T$, then 
$$ 
N (T) =  \frac{T}{2\pi} \log {\frac{T}{2\pi} } - \frac{T}{2\pi} + O(\log T)
$$ 
as $T \to \infty$\,. 
\end{enumerate}

Since the zeros of $\zeta(s)$ are symmetric about the critical line, the 
simplest possible assumption is that they all lie on the line. This is the 
famous  \\

\noindent{\bf{Riemann Hypothesis:}} If $\rho = \beta + i \gamma$  is a nontrivial zero 
of the zeta--function, then $\beta = 1/2$.  \newline

I will discuss the evidence for the truth of the Riemann Hypothesis later.
First, however, I want to explain the most direct connection between the 
primes and zeros of the zeta--function.

\noindent \section*{2 \quad The Prime Number Theorem} 

The Prime Number Theorem is the fundamental statistical fact about the primes.\\

\noindent{\bf{Prime Number Theorem.}} Let $\pi (x) = \sum_{p \leq x} 1 $.
Then we have
$$
\pi (x) \sim \frac x{\log x} \qquad \hbox{as} \quad x \to \infty     \,.
$$ 

One interpretation of the theorem is that the probability that a positive integer 
chosen at random in the interval $[1, x]$ is a prime equals $1/\log x$. 
Another is that the average distance between consecutive primes in 
the interval $[1, x]$ is $\log x$. 
 
For technical reasons, it is more convenient to express the  theorem 
in the following form, which can be shown to be equivalent by partial 
summation. \\

\noindent{\bf{Prime Number Theorem (second version).}} Set 
 $$ 
 \Lambda (x) = \begin{cases} \log p & \hbox{ if } x =
 p^k, \\ 0 & \hbox{ if } x \ne p^k \;. \end{cases}
 $$
 Then we have
 $$
 \psi (x) = \sum_{n \leq x} \Lambda (n) \sim x  \qquad \hbox{as} \quad x \to \infty     \,.
 $$\\

The proof I'll sketch  here is based on the ``explicit formula'', 
which is called that because it explicitly shows the relationship between 
the zeros and primes.

We begin by assuming that  $\Re s > 1$.  From the Euler product 
representation for the zeta--function we see that
$$ 
 \log \zeta (s)  = - \sum_p \log (1-p^{-s}) = \sum_p \sum^\infty_{k=1}
 \frac{1}{k p^{ks}} . 
$$
Differentiating, we find that
$$ 
- \frac{\zeta'}\zeta (s)  =  
\sum_p \sum^\infty_{k=1} \frac{\log p}{p^{ks}} = \sum_n \frac{\Lambda (n)}{n^s}
 \;.
$$
Here we have used a consequence of the fact that 
$\zeta(s)$ has an Euler product, namely, that its logarithm and, 
therefore, its logarithmic derivative also have Dirichlet series representations.
 
The idea now is to express the sum up to $x$ of the coefficients
$\Lambda(n)$ of the last series  (that is, $\psi (x)$) as an integral transform.  
This is analogous to writing the Fourier coefficients of a 
periodic function as an integral.

We break the argument into steps. \\ 
\noindent{\bf{StepI.}} Note that
$$\frac 1{2\pi i} \int^{2 + i\infty}_{2-i\infty} \frac{y^s}s ds
= \begin{cases}
1 & \hbox{ if } y > 1 \;, \\
1/2 & \hbox{ if } y = 1\;, \\
0 & \hbox{ if } 0 < y < 1 \; . 
\end{cases}
$$
This is a standard exercise in complex function theory. If $y>1$ we 
may pull the contour left to $-\infty$. In doing so we pass a simple 
pole of the integrand with residue $1$. If $y<1$ we 
pull the contour right to $+\infty$. This time we pass no 
poles, so the value of the integral is  $0$. When $y=1$, we can 
calculate the Cauchy principal value of the integral directly, and it 
turns out to be $1/2$.

\noindent{\bf{StepII.}}\\
We use the formula above to evaluate
 $$
 \begin{aligned}
 \frac 1{2\pi i } \int^{2 + i \infty}_{2 - i \infty} \left( -
 \frac{\zeta'}\zeta (s) \right) \frac{x^s}s ds & = \frac 1{2\pi i}
 \int^{2 + i \infty}_{2 - i \infty} \left( \sum^\infty_{n=2}
 \frac{\Lambda (n)}{n^s} \right) \frac{x^s}s ds \\
 & = \sum^\infty_{n=2} \Lambda (n) \left( \frac 1{2\pi i } \int^{2 +
 i\infty} _{2 - i \infty} \; \frac{(x/n)^s}{s} ds \right) \\
 & = \psi (x) - \frac 12 \Lambda (x)  \,. \end{aligned}
 $$
The interchange of summation and integration is not quite justified 
here. We should really truncate the integral first and keep track of 
the error terms. But we will ignore this technical point so as not to 
obscure the main idea.

\noindent{\bf{StepIII.}} \\
Evaluate the integral in Step II in a different way by pulling the contour 
left to $ -\infty$. We pick up residues from the simple poles of 
$-\frac{\zeta'}\zeta (s) \frac{x^s}s$ at i) the trivial and nontrivial zeros of $\zeta (s)$, 
ii) the pole of $\zeta (s)$ at $s =1$, and iii) the pole at $s =0$. 
Calculating and summing the residues, and then equating the result to 
$\psi (x) - \frac 12 \Lambda (x)\,,$ we find that 
$$ 
\psi (x) - \frac 12 \Lambda (x) = x - \sum_\rho \frac{x^\rho}\rho -
\frac{\zeta'}\zeta (0) - \sum^\infty_{m=1} \frac{x^{- 2m}}{-2m} \,.
$$

This is the ``explicit formula''. Had we worked with a truncated integral over 
the interval $[2-iT, 2+iT]$, say, rather than the integral over $[2-i\infty, 
2+i\infty ]$ (as we should have done to overcome the convergence problem in Step 
II), the sum over $\rho$'s would also be truncated. The analysis is more 
complicated, but leads to a more useful form of the explicit 
formula, namely,
$$ 
\psi (x) = x - \sum_{|\gamma | \leq T} \frac{x^\rho}\rho +
\mathcal{E} (x, T) \,,
$$
where $\mathcal{E} (x, T)$ is a known error term. For the 
applications we have in mind, one can show that it is possible to 
choose $T$ as a function of $x$ in such a way that the error  
term is not significant. We therefore will not bother with the exact 
form of $\mathcal{E} (x, T)$.
 
From the last form of the explicit formula one can almost see
the Prime Number Theorem.  Since $|x^\rho|= |x^{\beta+i\gamma}|
= x^\beta$, the term involving the sum over zeros should be $o(x)$ as long as the $\beta$ are 
not too close to 1.  Indeed, using the estimate $N (T) <<  T \log T $,
we see that the sum is 
$$
<< (\max_{\substack{0< \gamma \leq T}} 
x^\beta ) \sum_{0< \gamma \leq T} |\rho|^{-1}
<< (\max_{\substack{0< \gamma \leq T}} x^\beta )\log^{2}T .
$$
Now, one can show that the inequality $\beta \leq 1 - 
\frac{c}{\log \gamma}$ holds, where $c$ is a positive constant, for every zero 
$\rho = \beta + i \gamma$.
This leads to the Prime Number Theorem with an error term:
$$
\psi (x) = x + O \left( xe^{-b\sqrt{\log x} } \right)\,,
$$ 
with $b$ a positive constant.
Clearly the farther left the zeros all lie from the line $\Re s=1$, the better the error
term. Since the zeros are symmetric about the line $\sigma = 1/2$, 
the farthest left they can be is on the critical line, and 
in this case one can show that the sum over zeros in the explicit 
formula is $O(x^{1/2+\epsilon})$. Thus the Riemann Hypothesis implies that 
$$
\psi (x) = x + O(x^{1/2+\epsilon}) \,.
$$
In fact, this statement also implies the Riemann Hypothesis.

Why do we care about the error term? Because the main term just tells
us the large scale behavior of the sequence of primes; all the detailed
fluctuations in the counting function for the primes is hidden in the $O$-term.  
To illustrate this point, let us assume RH and consider the problem of
how large the gaps between consecutive primes can be. From 
$\psi (x) = x + O (x^{1/2 + \epsilon})$ we easily see that
$$ 
\psi (x + h) - \psi (x) = h + O (x^{1/2 + \epsilon}), \quad 1 \leq h
\leq x\;.
$$
Now suppose that there are no primes in $[x, x + h)$ (it can easily be shown 
that we may ignore the prime powers). Then
$0 = h + O (x^{1/2 + \epsilon})$, so we have
$ h = O (x^{1/2 + \epsilon})\;.$
Thus, on RH there is a positive constant $C$ such that the interval
$(x, x + Cx^{1/2 + \epsilon}]$ always contains a prime.
Hence, the error term in the Prime Number Theorem has a bearing on the size of the 
maximal gap between primes. Had we not assumed RH, the same analysis would 
have only led to the 
assertion that every interval $(x, x + Cxe^{-b\sqrt{\log x} } ]$
contains a prime. This of course is much weaker.

\noindent \section*{3 \quad The evidence for the Riemann Hypothesis} 
 
I will conclude by indicating why we believe the  Riemann Hypothesis.
The main evidence supporting it is the following.

\begin{enumerate}
 
\item {\bf{ Zero--free regions:}} There is a region to the left of the line 
$\Re s=1$ that is free of zeros. More specifically, there is a positive 
constant $c$ such that the region in the critical strip bounded by 
the curve $\sigma =1-c/\log (|t|+2)$ on 
the left, and $\sigma =1$ on the right, contains no zero of $\zeta(s)$. 
(We used this fact when we deduced the Prime Number Theorem with error term.)
The region has been widened slightly, but no one has been able to 
extend it to a vertical strip.
The conjecture that there is such a strip is refered to as the Quasi--Riemann 
Hypothesis.

\item {\bf{ Zero--density estimates:}}
Let $N(\sigma, T)$ denote the number of zeros $\rho = \beta + i \gamma$ of the zeta--function 
such that $\sigma \leq \beta \leq 1$ and  $0 < \gamma \leq T$. 
Many estimates have been proved of the type $N(\sigma, T) \leq 
T^{\lambda (\sigma)}$ with 
$\lambda (\sigma) < 1$ and $\lambda (\sigma)$ decreasing for $1/2 < \sigma \leq 1$.

\item {\bf{ Calculations of zeros:}}
The first fifty billion zeros of the zeta--function above the real 
axis have been shown to be simple and to lie on the critical line. 
Also, A. M. Odlyzko [O] has  performed extensive computations showing, among many other 
things, that the nearest several hundred million zeros to the $10^{20}$th zero
lie on the critical line.
Zeros of many other $L$--functions have also been computed and 
all of these have been shown to lie on the (corresponding) critical line.

\item {\bf{ Estimates of zeros on the critical line: }}  
Let $N_{0}(T)$ denote the number of zeros of $\zeta(s)$ on the 
critical line whose ordinates $\gamma$ satisfy $0< \gamma \leq T$ .
In 1914, Hardy [H] showed that $N_{0}(T) \to \infty$ with $T$. 
In 1921, Hardy and Littlewood [HL2] showed that $N_{0}(T) > T$.
Then, in 1942, A. Selberg [S] proved that $N_{0}(T) \geq c N(T)$,
for some positive constant $c$.
Thus, a positive proportion of the zeros lie on the critical line.
The constant $c$ was quite small, but
in 1974 N. Levinson [L], using a different method, showed that $N_{0}(T) > 1/3  N(T)$. 
In 1989, B. Conrey [C], iincreased the proportion to more than $2/5$.

\item {\bf{ The finite field case:}}
It is possible to define analogues of the zeta--function for curves 
and varieties over finite fields. It has been shown that the 
analogous Riemann Hypotheses for these
zeta-functions are true.
\end{enumerate}
\newpage

%%%%%%%%
%%%%%%%%
\centerline{\bf\LARGE{Lecture II}} 
\medskip
 
\centerline{\bf{\Large{Mean-Value Theorems and the Zeros: 
Three Applications}}}
\medskip

\noindent \section*{1 \quad An introduction to mean-value formulas}

In this lecture I will explain what mean--value estimates are, 
give a sampling of some of the most important ones, and 
present three applications to the study of the zeros of the 
zeta--function. These should make it clear why they play such a 
central role in the theory.

Let's begin with some general remarks on mean--value theorems. 

By a mean--value theorem we mean an estimate for an integral of the 
type
$$
\int_{0}^{T} {| F(\sigma + i t) |}^{2}\,dt 
$$
or
$$
\int_{0}^{T}  F(\sigma + i t) \,dt 
$$
as $T \to \infty$, where $F(s)$ is a function representable by a 
convergent Dirichlet series in some half--plane $\Re s > \sigma_{0}$ of the 
complex plane. The path of integration here need \emph{not} lie in this 
half--plane. For example, we would like to know the size of the 
integrals 
$$
I_{k}(\sigma, T) = \int_{0}^{T} {| \zeta(\sigma + i t) |}^{2k}\,dt \,,
$$
for $\sigma \geq 1/2$ and $k$ a positive integer. Here 
$F(s) = \zeta(s)^{k}$ and its Dirichlet series converges only for 
$\sigma > 1$.

There are many variations on this theme. For example, one 
might also consider a discrete version, namely an estimate for a sum of 
the form 
$$
\sum_{r=1}^{R} {| F(\sigma_{r} + i t_{r}) |}^{2}\,,
$$
where the points $\sigma_{r} + i t_{r}$ lie in $\CC$.
Another possibility is for $F(s)$ to involve a parameter $N$, say.
We then desire as uniform an estimate as possible in 
both $N$ and $T$. The simplest case is when 
$$
F(s) = F_{N}(s) = \sum_{n=1}^{N} a_{n} n^{-s}
$$
is a Dirichlet polynomial. Here one can calculate the mean--value 
in a straightforward way. We have
\begin{equation*}
\begin{aligned}
\int_{0}^{T} {| F_{N}(\sigma + i t) |}^{2}\,dt
& = \int_{0}^{T} |\sum_{n=1}^{N} a_{n} n^{-\sigma - i t}|^{2}dt \\
& = \int_{0}^{T}\sum_{n=1}^{N} \sum_{m=1}^{N} a_{n} \overline{a}_{m}
  n^{-\sigma - i t}m^{-\sigma + i t}\,dt \\
& = \sum_{n=1}^{N} \sum_{m=1}^{N} \frac{a_{n} \overline{a}_{m}}
{(n m)^{\sigma}} \int_{0}^{T} (m/n)^{it}\,dt  \\
& = T \sum_{n=1}^{N} \frac{|a_{n}|^{2}}{n^{2\sigma}} +
\sum_{\substack{1\leq m, n \leq N \\ m\neq n}}
\frac{a_{n} \overline{a}_{m}}
{(n m)^{\sigma}} (\frac{ (n/m)^{iT} - 1}{i \log n/m}) \,.
\end{aligned}
\end{equation*}
It is not difficult to show that the second term on the last line 
is
$$
\ll N\log N \sum_{n=1}^{N} \frac{|a_{n}|^{2}}{n^{2\sigma}}.
$$

Hence, we find that
$$
\int_{0}^{T} |\sum_{n=1}^{N} a_{n} n^{-\sigma - i t}|^{2}dt 
= (T + O(N\log N)) \sum_{n=1}^{N} \frac{|a_{n}|^{2}}{n^{2\sigma}} \,.
$$
From this we see that, as long as $N \ll T^{1-\epsilon}$ for some 
small positive $\epsilon $, the asymptotic estimate
$$
\int_{0}^{T} |\sum_{n=1}^{N} a_{n} n^{-\sigma - i t}|^{2}dt 
= T  \sum_{n=1}^{N} \frac{|a_{n}|^{2}}{n^{2\sigma}} 
$$
holds. On the other hand, when $N \gg T$ the mean-value can be about 
as large as 
$$
N\log N \sum_{n=1}^{N} \frac{|a_{n}|^{2}}{n^{2\sigma}}.
$$
Thus, the size of the mean-value is dominated by the contribution 
of ``diagonal''  terms when $N$ is smaller than $T$ but, in the 
opposite case, the main contribution may be from the ``off--diagonal'' 
terms. Goldston and Gonek [GG] have given a much more precise version of 
the mean--value formula for such ``long'' Dirichlet polynomials 
in terms of the size of the coefficient 
correlations sums $\sum_{n=1}^{N} a_{n} a_{n+h}$.

\noindent \section*{2 \quad Connections between zeros and mean-values}

Mean--value estimates are used in many ways to study the zeros of the 
zeta--function; indeed, this is one of the reasons that so much effort has been 
expended on them.  Why should there be a connection? One direct link 
is the general relationship between the zeros of an analytic function and its average 
size as given by Jensen's Formula in classical function theory.

\paragraph{{Jensen's Formula:}} Let $f(z)$ be analytic for $|z|\leq R$
and suppose that $f(0)\neq 0$. If $r_{1}, r_{2}, \ldots, r_{n}$ are the moduli of 
all the zeros of $f(z)$ inside $|z|\leq R$, then
$$
\log(\frac{|f(0)| R^{n}}{r_{1}r_{2}\cdots r_{n}}) =
\frac{1}{2\pi}\int_{0}^{2\pi}\log |f(R e^{i\theta})|\,d\theta \,.
$$
Here we see that the size of the mean--value of $\log|f(z)|$, 
this time around a circle, is related to the distribution of 
the zeros of $f(z)$ inside that circle. There is an analogous result 
for rectangles, which is often  more useful when working with 
Dirichlet series, namely,

\paragraph{{Littlewood's Lemma:}} Let $f(s)$ be analytic and nonzero 
on the rectangle $\mathcal{C}$ with vertices $\sigma_{0}$, 
$\sigma_{1}$, $\sigma_{1}+iT$, and $\sigma_{0}+iT$, where 
$\sigma_{0} < \sigma_{1}$.
Then
\begin{equation*}
\begin{aligned}
2\pi \sum_{\rho \in \mathcal{C}} \hbox{Dist}(\rho) 
&= \int_{0}^{T} \log 
|f(\sigma_{0} +it)|\,dt - 
\int_{0}^{T} \log |f(\sigma_{1} +it)|\,dt \\
& +\int_{\sigma_{0}}^{\sigma_{1}} \arg 
f(\sigma+iT)\,d \sigma 
-\int_{\sigma_{0}}^{\sigma_{1}} \arg 
f(\sigma)\,d\sigma \,,
\end{aligned}
\end{equation*}
where the sum runs over the zeros  $\rho$ of $f(s)$ in $\mathcal{C}$
and  ``$\hbox{Dist}(\rho)$'' is the distance from $\rho$ to the left 
edge of the rectangle.

When we use Littlewood's Lemma below, it will turn out that only the 
first term on the right--hand side is significant. So in order to not 
get too technical, I will always use the result in the form
$$
2\pi \sum_{\rho \in \mathcal{C}} \hbox{Dist}(\rho) 
= \int_{0}^{T} \log 
|f(\sigma_{0} +it)|\,dt  + \mathcal{E} ,
$$

where $\mathcal{E}$ is an error term that can be ignored and may be different 
on different occassions.

The Integral of the logarithm usually cannot be dealt with directly, 
so we often use the following trick.  We have
$$
\frac{1}{T} \int_{0}^{T} \log |f(\sigma + it)|\,dt 
=\frac{1}{2T} \int_{0}^{T} \log (|f(\sigma + it)|^{2}) \,dt 
\leq \frac{1}{2} \log (\int_{0}^{T} |f(\sigma + it)|^{2} \,dt ) \,,
$$
where the inequality follows from the arithmetic--geometric mean 
inequality. In this way we see a direct connection between the location of 
the zeros within a rectangle and the type of mean--values we have 
been considering.

\noindent \section*{3 \quad A sampling of mean--value results}
A great deal of work has been devoted to estimating the means
$$
I_{k}(\sigma, T) = \int_{0}^{T} {| \zeta(\sigma + i t) |}^{2k}\,dt \,.
$$
When $k=1$ we know that for each fixed $\sigma >1/2$ 
$$
I_{1}(\sigma, T) \sim c(\sigma)\, T \,,
$$
as $T \to \infty$, where $c(\sigma)$ is a know function of $\sigma$.
In 1918 Hardy and Littlewood [HL1] proved that, on the critical line itself, 
$$
I_{1}(1/2, T) \sim T\log T \,.
$$
What can such estimates tell us about the zeta--function?
Comparing the result for $\sigma$ greater than $1/2$ with that for
$\sigma = 1/2$, we see that the zeta--function tends to assume, on 
average, much larger values on the critical line than to the right 
of it. Since it 
also has many zeros on the critical line, we should expect the 
zeta--function to behave rather erratically there. 

The next higher moment was detrermined in 1926 by Ingham [I1], who proved that 
$$
I_{2}(1/2, T) \sim \frac{T}{2{\pi}^{2}}\log^{4} T \,.
$$
Unfortunately, no asymptotic estimate for any $k$ greater than $2$ has ever been 
proven.  Ramachandra [R] has shown that 
$$
I_{k}(1/2, T) \gg T {\log}^{k^{2}} T \,,
$$
and we expect that 
$$
I_{k}(1/2, T) \ll T \log^{k^{2}} T \,.
$$
Conrey and Ghosh [CG1] have conjectured that 
$$
I_{k}(1/2, T) \sim \frac{a_{k}g_{k}}{\Gamma(k^{2}+1)} \, T \log^{k^{2}} T \,,
$$
where 
$$
a_k= \prod_p \left( \left( 1 - \frac 1p
\right)^{(k-1)^2} \sum^{k-1}_{r=0} \binom{k-1}{r}^2 p^{-r} \right)
$$ 
and 
$g_{k}$ is an unknown constant. Not only does a proof of the 
conjecture seem far off, but it is only recently that anyone been able 
to suggest a plausible value for $g_{k}$. I will return to the 
problem of $g_{k}$ in the final lecture.

Another type of mean--value important in applications is 
$$
\int_{0}^{T} { |\zeta(\sigma + i t) M_{N}(\sigma + i t)| }^{2}\,dt \,,
$$
where 
$$
M_{N}(s) = \sum_{1\leq n \leq N} \frac{\mu(n)}{n^{s}}  
P(\frac{\log n}{\log N}) 
$$
and $P(x)$ is a polynomial.
Since 
$$
\frac{1}{\zeta(s)} = \sum_{n=1}^{\infty} \frac{\mu(n)}{n^{s}}  
\qquad(\Re s > 1) \,,
$$
we can view $M_{N}(s) $ as an approximation to the reciprocal of 
$\zeta(s)$ in $\Re s > 1$. We might then expect the approximation
to hold (in some sense) inside the critical strip as well. If that is the case, we 
should also expect that multiplying the zeta--function by $M_{N}$ dampens, 
or mollifies, the large values of zeta. Below we will see two
applications of this idea. The  most general estimates known for such
integrals are due to Conrey, Ghosh, and Gonek [CGG2], who obtained asymptotic 
estimates for them when the length of the Dirichlet 
polynomial $M_{N}(s)$ is $N=T^{\theta}$ with $\theta < 1/2$. Later, 
Conrey [C] used Kloosterman sum techniques to show that these formulas 
also hold for $\theta < 4/7$.

Assuming the Riemann Hypothesis and the Generalized Lindeloff 
Hypothesis are true, Conrey, Ghosh, and Gonek [CGG2] 
also proved discrete versions of such 
mollified mean--values, including estimates for sums of the type
$$
\sum_{0< \gamma <T} {|{\zeta}^{'}(\rho) M_{N}(\rho)| }^{2} \,,
$$ 
where $\gamma$ runs over the ordinates of the zeros 
$\rho = \beta + i\gamma $
of $\zeta(s)$. The first result of this type, but without 
the polynomial $M_{N}$, were proved by Gonek[G] under the 
assumption of the Riemann Hypothesis alone.

Having presented a brief catalogue of mean--value estimates, I will 
now turn to a few of their applications.

\noindent \section*{4 \quad A simple zero--density estimate}

We want to show that there are relatively few zeros of the 
zeta--function in the right half of the 
critical strip. Let $\sigma_{0}$ be a fixed real number strictly between $1/2$ 
and $1$ and let $\mathcal{C}$ be the rectangle in the complex plane with 
vertices at $2$, $2+iT$, $\sigma_{0} +iT$, $\sigma_{0}$. 
Applying our (simplified) version of Littlewood's Lemma, we see that
$$
\sum_{\rho \in \mathcal{C}} \hbox{Dist}(\rho)  = \frac{1}{2\pi} \int_{0}^{T} \log 
(|\zeta(\sigma_{0}+it)|) \,dt + \mathcal{E} ,
$$
Where  $\hbox{Dist}(\rho)$ is the distance of the zero $\rho$ from the line 
$\Re s=\sigma_{0}$. Now let $\sigma$ be a fixed real number with 
$\sigma_{0} < \sigma <1$ and write $N(\sigma, T)$ for the number of 
zeros $\rho = \beta + i\gamma$ of $\zeta(s)$ with $\sigma < \beta \leq 2$ and 
$0<\gamma< T$. On the one hand, we have 
$$
\sum_{\rho \in \mathcal{C}}  \hbox{Dist}(\rho)  \geq \sum_{\substack{\rho \in 
\mathcal{C}\\ \sigma \leq \beta }} \hbox{Dist}(\rho) 
\geq (\sigma - \sigma_{0}) N(\sigma, T).
$$
On the other hand, 
\begin{equation*}
\begin{aligned}
\frac{1}{2\pi} \int_{0}^{T} \log 
(|\zeta(\sigma_{0}+it)|) \,dt
&=
\frac{1}{4\pi} \int_{0}^{T} \log 
(|\zeta(\sigma_{0}+it)|^{2}) \,dt \\
&\leq
\frac{T}{4\pi} \log (\frac{1}{T}
\int_{0}^{T} |\zeta(\sigma_{0}+it)|^{2}) \,dt
\end{aligned}
\end{equation*}
by the arithmetic--geometric mean inequality, as 
before. The integral on the last line is $I_{k}(\sigma_{0}, T)$, which 
we have seen is $\sim c(\sigma_{0}) T$, where $c(\sigma_{0})$ is 
positive and independent of $T$. 
Thus, the last expression is $O(T)$ . It follows that 
$$
N(\sigma, T) \ll T \,.
$$
Since $N(T) \sim \frac{T}{2\pi}\log T$, we 
see that
$$
N(\sigma, T)/N(T) = O(\frac{1}{\log T})
$$
for any fixed $\sigma > 1/2$. We may interpret this as saying
that the proportion of zeros to the right of any line 
$\Re s = \sigma > 1/2$ is infinitesimal. 

This, the first zero--density estimate, was
proved by H. Bohr and E. Landau [BL] in 1914. Since then 
much stronger results have been proven, typically  
of the form
$$
N(\sigma, T) << T^{ \lambda(\sigma)}\,,
$$
where $\lambda(\sigma) < 1$ for $\sigma > 1/2$.
Nevertheless, the underlying idea in the proof of many (but not all) of these 
results already appears here. 
 
\noindent \section*{5 \quad Levinson's method}

Zero--density theorems tell us there are (relatively) few zeros to the right of the 
critical line. Our goal here is to sketch the metod of Levinson [L], which
shows that there are many zeros \emph{on} it. 

Recall that
\begin{equation}
\begin{aligned}
N(T) & = \# \{ \rho = \beta + i\gamma \mid \zeta (\rho) =0, \quad 0 <
\gamma < T \} \notag\\
& \sim \frac T{2\pi} \log T \notag\\
\end{aligned}
\end{equation}
and let
$$
N_0 (T)  =  \# \left\{ \rho = \frac 12 + i\gamma \mid \zeta (\rho ) = 0,
\quad 0 < \gamma < T \right\} 
$$
denote the number of zeros  on the critical line up to height $T$. 
The  important estimations of $N_0 (T)$ were: \newline
G. H. Hardy (1914) : \quad $N_0 (T) \to \infty$ \qquad (as $T\to
\infty$) \newline
G. H. Hardy-J. E. Littlewood (1921) : \quad $N_0 (T) > cT$ \newline
A. Selberg (1942) : \quad $N_0 (T) > c' N(T)$ \newline
N. Levinson (1974) : \quad $N_0 (T) > \frac 13 N(T)$ \newline
J. B. Conrey (1989) : \quad $N_0 (T) > \frac 25 N(T)$ \newline

In keeping with the  theme of this lecture, I should point out that each of 
the last four results requires the use of mean--value theorems.

Levinson's method begins with the following fact first proved by 
Speiser [Sp].

\noindent{\bf{Theorem (Speiser).}} The Riemann Hypothesis is equivalent 
to the assertion that $\zeta^{\prime}(s)$ does not vanish in the left half 
of the critical strip.

In the early seventies, N. Levinson and H. L. Montgomery [LM]  proved a 
quantitative version of this. Let 
$$
N_{-}^{\prime}(T) = \# \left\{ \rho^{\prime} = \beta^{\prime} + 
i\gamma^{\prime} \mid \zeta^{\prime} (\rho ^{\prime} ) = 0, \,
-1 < \beta^{\prime} < 1/2, \quad 0 < \gamma ^{\prime}  < T \right\} 
$$
and 
$$
N_{-}(T) = \# \left\{ \rho = \beta + 
i\gamma \mid \zeta(\rho ) = 0, \,
-1 < \beta < 1/2, \quad 0 < \gamma < T \right\} \,.
$$

\noindent{\bf{Theorem (Levinson-Montgomery).}}  We have 
$ N_{-}(T) = N_{-}^{\prime}(T) + O(\log T)$ .

The idea behind the proof is as follows. Let $0<a<1/2$ and let 
$\mathcal{C}$ denote the positively oriented rectangle with vertices
$a+iT/2$,  $a+iT$,  $-1+iT$, and $-1+iT/2$ . By a standard method it 
is not difficult to show that 

$$ 
\Delta \hbox{ arg } \frac{\zeta'}{\zeta} (s) \bigg|_{\mathcal{C}}
= O (\log T),
$$
independently of $a$.
Given this, we see that 
$$ 2 \pi( \# \hbox{ zeros of } \zeta' (s)\,  
\hbox{in}\, \mathcal{C} - \# \hbox{ zeros of } \zeta
(s) \,\hbox{in} \,\mathcal{C})
= O(\log T) .
$$ 
The theorem now follows on observing that $a$ was arbitrary,
and by ``adding'' rectangles with top and bottom edges, respectively, at
$T$ and  $T/2$, $T/2$ and  $T/4$, \ldots .

We now sketch Levinson's method. We have just seen that 
$ N_{-}(T) = N_{-}^{\prime}(T) + O(\log T)$.  Now, the nontrivial zeros of 
$\zeta(s)$ are symmetric about the critical line. Hence, the  number of 
them lying to the right of the critical line, to the left of the line 
$\sigma=2$, and above the real axis up to 
height $T$ is also
$N_{-}(T)$. Therefore

\begin{equation*}
\begin{aligned}
N(T) &= N_{0}(T) + 2 N_{-}(T) \\
&= N_{0}(T) + 2 N'_{-}(T) + O (\log T) ,
\end{aligned}
\end{equation*}
or
$$
 N_0 (T) = N(T) - 2 N'_{-}(T) + O (\log T) .
$$

The size of the first term on the left hand side of the last 
line is known, namely, \\$(1+o(1)) \frac T{2\pi} \log T$. Hence, if we 
can determine a sufficiently small upper bound for  $N'_{-}(T)$, we 
can deduce a lower bound for $ N_0 (T)$. 

To find such an upper bound it is convenient to first note that 
the zeros of $\zeta' (s)$  in the region $-1 < \sigma <
1/2, \; 0 < t < T$,  are identical to the zeros of
$\zeta' (1-s)$ in the reflected region $1/2 < \sigma < 2, \; 0 <
t < T$ . One can also show, by the functional equation
of the zeta--function, that $ \zeta' (1-s)$ and $G(s) = 
\zeta (s) + \zeta' (s) / L(s)$, where $L(s)$ is essentially 
$\frac{1}{2\pi} \log T$, have the same zeros in  $1/2 < \sigma < 2, \; 0 <
t < T$ . It turns out to be technically advantageous to count the 
zeros of $G(s)$ rather than those of $\zeta'(1-s)$. 

To bound the number of zeros of $G(s)$ in this region, we apply Littlewood's Lemma. 
Let $a = \frac 12 - \frac\delta{\log T}$, with $\delta$ a small 
positive number, and let $\mathcal{R}_{a}$ denote the rectangle whose 
vertices are at $a, 2, 2+iT$, and $a+iT$. It would be natural to apply 
our abreviated form of the lemma to obtain

$$\underset{\rho^{*} \in \mathcal{R}_{a}}{\sum \;\; }
\hbox{Dist}(\rho^{*}) 
= \frac 1{2\pi} \int^T_0 \log | G(a+ it) | dt + \mathcal{E},
$$
where $\rho^{*}$ denotes a zero of $G(s)$ and Dist($\rho^{*}$)
is its distance to the left edge of $\mathcal{R}_{a}$. However, in 
the next step, when we apply the arithmetic--geometric mean inequality to 
the integral, we would lose too much. To avoid this loss, we first dampen, or 
mollify, $G(s)$ and apply Littlewood's Lemma in the form
$$
\sum_{\substack{\rho^{**} \in \mathcal{R}_{a}\\ GM(\rho^{**}) = 
0 }}
\hbox{Dist}(\rho^{**}) 
= \frac 1{2\pi} \int^T_0 \log | G(a+ it) M (a+it) | dt + \mathcal{E}\,.
$$
Here
$$ 
M(s) = \sum_{n\leq T^\theta} \; \frac{a_n}{n^s} \;,\qquad 
a_n = \mu(n)n^{a-1/2} \left( 1 - \frac{\log n}{\log T^\theta} \right) \,,
$$
approximates $1/\zeta(s)$ and $\theta > 0$. Note that included among the zeros of $G(s)M(s)$ 
in $\mathcal{R}_{a}$ are all the zeros of $G(s)$ in $\mathcal{R}_{a}$. 
Therefore we have
\begin{align*}
\sum_{\substack{\rho^{**} \in \mathcal{R}_{a}\\ GM(\rho^{**}) = 
0 }} \hbox{Dist}(\rho^{**}) 
\geq 
&\sum_{\substack{\rho^{*} \in \mathcal{R}_{a}\\ G(\rho^{*}) = 
0 }} \hbox{Dist}(\rho^{*}) \\
\geq &\sum_{\substack{\rho^{*} \in \mathcal{R}_{a},\Re \rho^{*} > 1/2 
\\ G(\rho^{*}) = 0 }} \hbox{Dist}(\rho^{*}) \\
\geq &(\frac12 - a) N'_{-}(T) \,.
\end{align*}
We now see that
\begin{align*}
(1/2 -a) N'(T) & \leq \frac 1{2\pi} \int^T_0 \log |GM (a + it)| dt + \mathcal{E} \notag\\
& = \frac 1{4\pi} \int^T_0 \log |GM (a + it)|^2 dt +\mathcal{E} \notag \\
& \leq \frac T{4\pi} \log \left( \frac 1T \int^T_0 | GM (a + it) |^2 dt
\right) + \mathcal{E}\,. \notag 
\end{align*}
Thus, we require an estimate for
$$
\int^T_0 | GM (a + it)|^2 dt\,.
$$
This is similar to a mean--value we saw in Section 3. 
Levinson was able prove an asymptotic estimate for this 
integral when $\theta = 1/2 - \epsilon$
with $\epsilon$ arbitrarily small. The resulting upper bound for 
$N'_{-}(T)$ then led to the lower bound
$$ 
N_0 (T) > \left( \frac 13 + o(1) \right) N(T) .
$$
Much later, Conrey was able to establish an asymptotic estimate when 
$\theta = 4/7 - \epsilon$, which led to
$$ 
 N_0 (T) > \left( \frac 25 + o(1) \right) N(T) .
$$
The form of the asymptotic estimate in both cases is the same as a function 
of $\theta$,  and D. Farmer [F] has given various heuristic arguments that 
suggest it should remain true even when one takes $\theta$ arbitrarily large. 
From Farmer's conjecture it follows that

$$
 N_0 (T) \sim N(T)\;.
$$

Before concluding this section, we remark that had we introduced a mollifier 
into our proof of the Bohr--Landau result in the previous section, we would 
have obtained a much stronger zero--density estimate.

\noindent \section*{6 \quad The number of simple zeros}

Our final application demonstrates the use of discrete mean--value theorems.

Let 
$$
N_{s}(T) = \# \{ \rho = \beta + i\gamma \mid \zeta (\rho) =0, \zeta' (\rho) 
\neq 0, \quad 0 < \gamma < T    \}
$$
denote the number of simple zeros of the zeta--function in the 
critical strip with ordinates between $0$ and $T$.
It is believed that all the nontrivial zeros are on the 
critical line and simple, in other words, that $N(T) =N_{0}(T) =N_{s}(T)$
for every $T > 0$. 
In 1973, H. Montgomery [Mo], used his pair correlation 
method to show that if the Riemann Hypothesis is true, then at least $2/3$ of the zeros are simple.
In other words, 
$$
N_{s}(T)/ N(T) > 2/3
$$
provided that $T$ is sufficiently large. We will present his argument in the 
third lecture. Now, however, we briefly describe a different method of Conrey, Ghosh, and 
Gonek [CGG1], which shows that on the stronger hypotheses of RH and the Generalized Lindeloff 
Hypothesis, one can replace the $2/3$ 
above by $19/27 = .703\ldots$. 

By the Cauchy--Schwarz inequality, we 
have
$$
\Big| \sum_{0< \gamma <T} \zeta^{'}(1/2+i\gamma) M_{N}(1/2+i\gamma) \Big|^{2} 
\leq 
\Big(\sum_{\substack{0< \gamma \leq T\\ 1/2+i\gamma \, \hbox{is simple}}} 
1 \Big)
\Big(\sum_{0< \gamma <T} {|{\zeta}^{'}(\rho) M_{N}(\rho)| }^{2} \Big) 
\,,
$$
where $M_{N}(s)$ is a Dirichlet polynomial of length $N$  with 
coefficients similar, but not identical, to those of $M(s)$ in the last section.  
Its purpose is also similar: to mollify $\zeta^{'}(1/2+i\gamma)$ so 
as to minimize the loss in applyng the Cauchy--Schwarz inequality.
If one assumes RH, the sum on the left--hand side is easy to compute 
and turns out to be $\sim \frac{19}{24} N(T) \log T$. The sum on the right--hand side is 
much more difficult to treat, but one can 
show that if RH and GLH are true, then it is $\sim \frac{57}{64} N(T) {\log}^{2} T$.
Inserting these estimates into the inequality above and solving for 
$N_{s}(T)$ leads to the result stated. An elaboration of the method 
leads to the conclusion that, on the same hypotheses, at least
$95.5\%$ of the zeros of $\zeta(s)$ are either simple or double.

\newpage 

\centerline{\bf\LARGE{Lecture III}} 
\medskip  
\centerline{\bf{\Large{Beyond the Riemann Hypothesis}}}
\medskip

\noindent \section*{1 \quad Gaps between primes again}

In the first lecture we saw  that the Prime Number 
Theorem with error term implies that
if $\psi(x+h) -\psi(x) = 0$, then there is a positive constant
$c_{1}$ such that $ h \ll x e^{-c_{1}\sqrt{\log x}} $. We also
saw that if the Riemann Hypothesis is true, then
$ h \ll x^{1/2+\epsilon} $ for any positive $\epsilon$.
The prime powers higher than the first contribute
at most $O(x^{1/2})$ to $\psi(x+h) -\psi(x)$, so
another way to phrase this is that the size of the gap between 
any two consecutive primes $p$ and $p^{\prime}$ is 
$O(p e^{-c_{1}\sqrt{\log p}}) $ unconditionally, and 
$O(p^{1/2+\epsilon})$ on RH. On the other hand, the Prime Number 
Theorem tells us that the size of the \emph{average} gap 
between $p$ and $p^{\prime}$ is $\sim \log p$.
This suggests that if the primes behave ``randomly'', then
$p^{\prime}-p \ll p^{\epsilon}$, and the numerical evidence
does indeed support this.  

Here we have a problem for which even the assumption of
the Riemann Hypothesis does not seem to give the right answer. 
The question I want to begin with here is: Why?

The  answer is not difficult to find.  Consider again the explicit formula 
$$ 
\psi (x) = x - \sum_{|\gamma | \leq T} \frac{x^\rho}\rho +
\mathcal{E} (x, T) \,,
$$
where $\mathcal{E} (x, T)$ is a known error term
and, from now on, we assume the Riemann Hypothesis. 
If we apply the formula with the arguments $x+h$ and $x$ and subtract, 
we obtain
$$ \begin{aligned}
\psi (x +h) - \psi (x) & = h - \sum_{|\gamma | \leq T} \; \frac{(x
+h)^\rho - x^\rho}\rho + \mathcal{E} (x+h, T) -  \mathcal{E} (x, T)       \\
& = h - \sum_{|\gamma | \leq T} \left( \int^{x+h}_x u^{\rho -1} du
\right) + \mathcal{E}^{\prime} (x, h, T) \\
& = h - \int_x ^{x +h} \left( \sum_{|\gamma | \leq T} u^{i \gamma}
\right) \frac{du}{\sqrt{u}} + \mathcal{E}^{\prime} (x, h, T)  \;.
\end{aligned}
$$
There is likely to be a lot of cancellation in the
sum in the integrand. However, when we estimated 
the error term in the Prime Number Theorem,  we lost it all by putting
absolute values around the individual terms.  Clearly this cancellation 
depends completely on the distribution of the sequence of ordinates 
$\gamma$. In other words, on the \emph{vertical} distribution of the zeros of the 
zeta--function.
 
This example is not unique; it often happens that the strength 
of the Riemann Hypothesis, or even of the Generalized Riemann Hypothesis, 
is not sufficient to establish what we think is the ultimate truth in important 
arithmetical questions. We also often find that we need to understand the 
vertical distribution of the zeros of the zeta--function 
and $L$--functions. 

\noindent \section*{2 \quad Pair correlation}

Prior to the early seventies, such an understanding seemed beyond 
reach. Then, in 1973, Hugh Montgomery [Mo] found a way to study the distribution 
of the differences between all pairs of ordinates of zeros of the Riemann zeta--function, 
assuming RH is true.  \\
 
\noindent{\bf{Montgomery's Theorem:}} Assume the Riemann Hypothesis. Set
$w (u) = \frac{4 }{4 + u^2}\;,$
and for $\alpha$ real and $T \geq 2$ write
$$
 F(\alpha) =F(\alpha, T) = (\frac{T}{2\pi} \log T)^{-1}
 \sum_{0 < \gamma,
\gamma' \leq T} \; w(\gamma - \gamma')
\; T^{i \alpha (\gamma - \gamma' )}\;,
$$
where $\gamma$ and $\gamma^{\prime}$ run over  ordinates of zeros 
of the Riemann zeta-function. Then $F(\alpha)$
is real and an even function of $\alpha$. Moreover, for any 
$\epsilon > 0$ we have
$$
 F(\alpha) = (1+o(1)) T^{-2\alpha} \log T + \alpha + o(1)\,  
   \qquad (\hbox{as} T \to \infty)   \\
$$
uniformly for $0 \leq \alpha \leq 1-\epsilon$. 

It was later observed that $F(\alpha)$ is nonnegative.

Integrating $F(\alpha)$ against a kernel $\hat{r}(\alpha)$, we see 
that
$$ \begin{aligned}
(\frac{T}{2\pi} \log T) \int_{-\infty}^{\infty} 
F(\alpha) \hat{r}(\alpha)\,d{\alpha} 
&= \int_{-\infty}^{\infty}  \sum_{0 < \gamma,
\gamma' \leq T} \; w(\gamma - \gamma')
\; T^{i \alpha (\gamma - \gamma' )} \hat{r}(\alpha)\,d{\alpha} \\
&= \sum_{0 < \gamma, \gamma' \leq T} \; w(\gamma - \gamma')
\int_{-\infty}^{\infty}  T^{i \alpha (\gamma - \gamma' )} 
\hat{r}(\alpha)\,d{\alpha}  \\
&= \sum_{0 < \gamma, \gamma' \leq T} \; 
r((\gamma - \gamma') \frac{\log T}{2 \pi} ) w(\gamma - \gamma') \,,
\end{aligned}
$$
where $r$ is the inverse Fourier transform of $\hat{r}$, that is,
$$
r(u) = \int_{-\infty}^{\infty} \hat{r}(\alpha)
e^{2\pi i \alpha u} \,d \alpha \,.
$$
Thus, the integral of $F(\alpha)$ against  a kernel $\hat{r}$ produces 
a sum involving the inverse transform $r$ evaluated at the differences 
of pairs of ordinates. Since Montgomery's Theorem is only valid in the range $-1 < 
\alpha < 1$, one can only use kernels $\hat{r}(\alpha)$ supported 
on $(-1, 1)$. For example, assuming RH and taking 
$ 
\hat{r} (\alpha) = \max \{ 0, \beta^{-1} (1 - | \alpha / \beta|) \}
$ 
with $0< \beta <1$, one obtains
$$
\sum_{0 < \gamma, \gamma' \leq T} \;
( \frac{\sin((\beta/2) (\gamma - \gamma') \log T)}
{(\beta/2) (\gamma - \gamma') \log T}  )^{2} \,w(\gamma - \gamma')
\sim (\frac{1}{\beta} + \frac{\beta}{3}) 
 \frac{T}{2\pi} \log T \,.
$$

Montgomery used this to obtain a lower bound for the number of simple 
zeros of the zeta-function as follows. First observe that
$$
\sum_{\substack{0 < \gamma ,\, \gamma' \leq T\\ \gamma = \gamma' }}\,1
\leq \sum_{0 < \gamma, \gamma' \leq T} \;
( \frac{ \sin((\beta/2) (\gamma - \gamma') \log T)}
{(\beta/2) (\gamma - \gamma') \log T}  )^{2} \,w(\gamma - \gamma') \,.
$$
Taking $\beta = 1 - \epsilon$ in this, we obtain
$$
\sum_{\substack{0 < \gamma ,\, \gamma' \leq T\\ \gamma = \gamma' }}\,1
\leq (\frac{4}{3} + o(1))  \frac{T}{2\pi} \log T \,.
$$
Now, if the zero $\frac{1}{2} +i \gamma$ has multiplicity $m(\gamma)$, 
then each $\gamma$ occurs 
$m(\gamma)$ times in the sum on the left. Thus, we have
$$
\sum_{\substack{0 < \gamma ,\, \gamma' \leq T\\ \gamma = \gamma'} }\,1
=\sum_{0 < \gamma < T} m(\gamma) 
$$
and therefore
$$
\sum_{0 < \gamma < T} m(\gamma)
\leq (\frac{4}{3} + o(1))  \frac{T}{2\pi} \log T \,.
$$
Finally, we easily see that
$$
\sum_{\substack{0 < \gamma \leq T\\ 
\frac{1}{2} + i \gamma \,\hbox{is simple}}}1
\geq \sum_{0 < \gamma < T} (2 - m_{\gamma})
\geq (2-\frac{4}{3} + o(1)) \frac{T}{2\pi} \log T \,.
$$
Hence, if the Riemann Hypothesis is true, then at least
two--thirds of the zeros are simple. Although this is not quite as 
strong a result as that obtained in Lecture II, 
namely $(19/27 + o(1)) \frac{T}{2\pi} \log T $, the 
hyotheses are also not as strong. For there we needed to assume
the Generalized Lindeloff Hypothesis in addition to RH.

Since we have focused so much on mean--value theorems, I should point 
out that Montgomery proved his theorem by relating $F(\alpha)$ to the  
mean--value of a Dirichlet series, namely,
$$
\frac{1}{x} \int_{0}^{T} |\sum_{n \leq x} \Lambda(n) 
(\frac{x}{n})^{-1/2 +it}   
+ \sum_{n> x} \Lambda(n) 
(\frac{x}{n})^{3/2 +it} |^{2} \,dt\,,
$$
where $x=T^{\alpha}$. Here we see a different explicit connection 
between the zeros and the primes. Indeed, Montgomery's starting point was 
a generalization of the explicit formula we saw in Lecture I (and 
again at the beginning of this lecture). 
The restriction $\alpha < 1$ in Montgomery's Theorem arises for a familiar reason: when 
$\alpha \geq 1$, the off--diagonal terms in the integral above contribute to the main term 
in the mean--value estimate. To determine this contribution (heuristically), Montgomery used 
a strong form of the Hardy--Littlewood twin prime conjecture. 
In this way he arrived at \\

\noindent{\bf{Montgomery's Conjecture:}} We have 
$$
 F(\alpha, T) = (1+o(1))    \qquad (\hbox{as} T \to \infty)
$$
for $\alpha \geq 1$, uniformly in bounded intervals.

This together with Montgomery's theorem determines $F(\alpha)$ on  all 
of $\mathbb{R}$.
Thus, one may use the conjecture to integrate $F(\alpha)$ against a much wider class of kernels 
than just those supported in $(-1, 1)$. Using an appropriate 
kernel he arrived at the \\

\noindent{\bf{Pair Correlation Conjecture:}}
For any fixed $\alpha$ and $\beta$
with  $ \alpha < \beta $, we have
$$
\sum_{\substack{0 < \gamma ,\, \gamma' \leq T\\ 2\pi\alpha/ \log T 
\leq \gamma' - \gamma \leq 2 \pi\beta / \log T} }1 \sim \left(
\int^\beta_\alpha 1 - \left( \frac{\sin \pi x}{\pi x} \right)^2 dx
+ \delta(\alpha, \beta )   \right) \frac T{2\pi} \log T
$$
as $T$ tends to infinity, where $\delta(\alpha, \beta ) =1$ 
if $0 \in [\alpha, \beta ]$, and  $\delta(\alpha, \beta ) =0$ 
otherwise. \\

The Pair Correlation Conjecture is an assertion about the distribution of the set 
of all differences between pairs of ordinates of the zeros. 
An enormous amount of data concerning the zeros has been collected and analyzed by 
A. M. Odlyzko [O], and the fit with the conjecture is remarkable.

As an example of the type of information we can deduce from it, 
let $0 < \alpha < \beta$ with $\alpha$ arbitrarily small. Then we find 
that
$$
\sum_{\substack{0 < \gamma ,\, \gamma' \leq T\\
 0 < \gamma' - \gamma \leq 2 \pi\beta / \log T} }1 \sim \left(
\int_{0}^{\beta} 1 - \left( \frac{\sin \pi x}{\pi x} \right)^2 dx
 \right) \frac T{2\pi} \log T \,.
$$
This shows that an infinite number of the zeros have another zero no 
farther away than $ 2 \pi\beta/ \log T$, no matter how small 
$\beta$ is.  We also deduce that 
$$
\sum_{\substack{0 < \gamma ,\, \gamma' \leq T\\ -2\pi\beta/ \log T 
\leq \gamma' - \gamma \leq 2 \pi \beta / \log T} }1 \sim \left(
\int^\beta_{-\beta} 1 - \left( \frac{\sin \pi x}{\pi x} \right)^2 dx
+ 1 \right) \frac T{2\pi} \log T\,.
$$
Combining this with the previous formula, we obtain
$$
\sum_{\substack{0 < \gamma ,\, \gamma' \leq T\\  
 \gamma' = \gamma } }1 \sim \frac T{2\pi} \log T\,.
$$
By our earlier discussion, we may write this as
$$
\sum_{0 < \gamma < T} m(\gamma)\sim \frac T{2\pi} \log T\,.
$$
On the other hand, von Mangoldt's formula tells us that 
$$
\sum_{0 < \gamma < T} 1 \sim \frac T{2\pi} \log T\,.
$$ 
It threfore follows that
$$
\sum_{\substack{0 < \gamma \leq T\\ 
\frac{1}{2} + i \gamma \,\hbox{is simple}}}1 \sim \frac T{2\pi} \log T\,.
$$
In other words, almost all the zeros are simple.

Before moving on we mention that
D. Goldston and H. Montgomery [GM] have shown that the Pair Correlation Conjecture 
is equivalent to a certain estimate of the variance of the number of 
primes numbers in short intervals.  D. Goldston, S. Gonek, and H. Montgomery 
[GGM] have shown that
it is also equivalent to an estimate for the mean--value
$$ 
\int^T_0 \left| \frac{\zeta'}\zeta ( \sigma + it) \right|^{2} \,dt \,,
$$
for $\sigma$ near $1/2$. Estimates of $F(\alpha, T)$ when $\alpha 
\geq 1$ remain elusive. The only progress in this 
direction so far is the lower bound  $F(\alpha, T) \geq 3/2 -\alpha + 
o(1)$ on the interval $ (1, 3/2)$ under the assumption of the Generalized Riemannn 
Hypothesis. This is due to D. Goldston, S. Gonek, A. E. \"{O}zl\"{u}k, and C. 
Snyder [GGOS].

\noindent \section*{3 \quad Random matrix theory}

Shortly after completing the work described above, Montgomery was told by F. Dyson 
that the ``form factor'' $1- \left( \frac{\sin \pi x}{\pi x} \right)^2$  
in the distribution law he had conjectured for pairs of zeros of the zeta--function 
is the same one that holds for pairs of 
eigenvalues of large random Hermitian matrices from 
the Gaussian Unitary Ensemble, or GUE, which we describe below. 
This and other matrix ensembles had been studied by physicists for decades because they 
can be used to model the Hamiltonians of complicated physical systems. 
The spectra, or energy levels, of such systems are given by the 
eigenvalues of the corresponding Hamiltonian. But in complicated 
situations, the Hamiltonian, let alone its eigenvalues, may not be known
with any certainty. In such cases the Hamiltonian can be modeled by  
large random Hermitian matrices with symmetry properties dictated by 
the physical situation. It is found that the average behavior of the 
eigenvalues of such families of matrices is often in agreement with 
the experimental data. Physicists are particularly interested in knowing 
various statistics of the energy levels, and pair correlation is 
merely one of these. They had also worked out ``$n$--level'' correlations of the 
eigenvalues, and Montgomery conjectured that the analogous law (there 
is a normalization one has to take into account) holds for the 
``$n$--level'' correlations of the zeros. Specifically, we have 
\\ 

\noindent{\bf{Montgomery's GUE Hypothesis:}} 
The distribution of
all $(n-1)$--tuples
$(\gamma_{2} -\gamma_{1}, \gamma_{3} -\gamma_{1}, \ldots , \gamma_{n} -\gamma_{1})$,
with the $\gamma_{i}$  ordinates of the zeros, has the form factor $\det K (x_1, \dots, x_n)$, where
$$
K(x_1, \dots, x_n) \; = \; (k_{i,j})^n_{i,j=1}\;, 
\qquad 
 k_{ii} = 1, \; \qquad  k_{ij} =
\frac{\sin \pi (x_i - x_j)}{\pi (x_i - x_j)} \;.\\
$$

The Pair Correlation Conjecture is the $n = 2$ case. 
Odlyzko [O] also used his data (alluded to above) to check 
this prediction, and the evidence is again compelling.
Moreover, so is the theoretical support 
(see, for example, E. Bogomolny 
and J. Keating [BK], D. Hejhal [He], and Z. Rudnik and P. Sarnak [RS]). 

Finally, the Gaussian Unitary Ensemble of order $N$
is the set of all $N \times N$ Hermitian matrices 
$H = (H_{j,k} )_{1\leq j,k\leq N}$ 
made into a probability space by equipping it with a probability 
measure $p(H) dH$, invariant under conjugation by all $N \times N$
Unitary matrices, where
$$
dH = \prod_{j \leq k}  d\Re H_{j,k}  \prod_{j < k} d\Im H_{j,k} \,
$$
and
$$
p(H) = \prod_{j \leq N}  \frac{1}{\sqrt{\pi}}e^{-H_{j,j}^{2}}    
\prod_{j < k}  \frac{2}{\pi}e^{-2((\Re H_{j,k})^{2}+ (\Im 
H_{j,k})^{2} )} \,.
$$

In practice it is often easier to work with the so called Circular 
Unitary Ensemble, or CUE, rather than the GUE. This is the
compact group of $N \times N$ unitary matrices equipped with Haar measure normalized 
so that the measure of the group is $1$. All eigenvalues have modulus one and the statistics of 
the eigenangles are known to be the same as those for 
the GUE eigenvalues. 

\noindent \section*{4 \quad Applications of random matrix theory to the zeta--function.}

Another remarkable development in the application of random matrix 
theory to analytic number theory has been the discovery by J. Keating 
and N. Snaith [KS] that the characteristic polynomial of a large random matrix 
from the Gaussian Unitary Ensemble or Circular Unitary Ensemble 
can be used to model the Riemann zeta--function and other L--functions. 

The idea is as follows.
Since Riemann's  function $\xi(s)$ is entire, it has a Hadamard product 
representation. Moreover, $\zeta(s)$ and $\xi(s)$ are the same up to 
well understood multiplicative factors. Therefore, one might plausibly
assume that at a large height $t$ in the critical strip, $\zeta(s)$ (with $s= 
\sigma + it$) should behave like a polynomial with the same zeros near $t$. 
If the zeros  are distributed like the 
eigenangles of matrices from the Circular Unitary Ensemble, one might 
then expect 
$$
Z_{N}(U,\theta) = \prod_{n=1}^{N} ( 1 - e^{i(\theta_{n} - \theta)} )\,,
$$
where the $\theta_{n}$ are the eigenangles of a random $N \times N$
unitary matrix $U$ from CUE, to model $\zeta(1/2 + it)$. For scaling reasons 
one takes $N = \frac{1}{2 \pi}\log t$. 

Keating and Snaith conjecture that 
the average of $|Z_{N}(U,\theta)|^{2k}$
over the full Circular Unitary Ensemble, with respect to Haar measure 
on the group, should be directly related to the $2k$th moment
$$
I_{k}(1/2, T)  = \int_{0}^{T}\, |\zeta(\frac{1}{2} + i t) |^{2k} \, dt
$$
of the zeta--function. Similarly,
the distribution of values of
$\log Z_{N}(U,\theta)$, say, should be the same as that of $\log \zeta(\frac{1}{2} + 
i t)$\,. The agreement with known results in both cases is 
remarkable. 

Consider the case of $I_{k}$. Recall from Lecture II that it had long been conjectured 
that there is a constant $ c_{k}$ such that
$$
I_{k}(1/2, T)  \sim  c_{k} T \log^{k^2} T
$$
as $T \to \infty$.  
J. B. Conrey and A. Ghosh [CG] have recast the conjecture into a more precise form,
namely that 
$$
c_{k} =  g_k a_k / \Gamma (k^2 +1) \,,
$$
where 
$$
a_k= \prod_p \left( \left( 1 - \frac 1p
\right)^{(k-1)^2} \sum^{k-1}_{r=0} \binom{k-1}{r}^2 p^{-r} \right)
$$ 
and $g_k$ is an integer. Thus the question comes down to the value of $g_k$.
The only proven values are the classical ones due to 
Hardy and Littlewood [HL1] and Ingham [I1] of
$g_1=1$ and $g_2=2$\,, respectively. Conrey and Ghosh [CG2] conjectured that 
$g_3 = 42$ and, using long Dirichlet polynomials to 
approximate $\zeta(s)^{k}$,  Conrey and Gonek [CGo] conjectured that $g_4 = 24024$.
At about the same time, Keating and Snaith [KS] calculated arbitrary complex 
moments of the characteristic polynomials $Z_{N}(U,\theta)$ averaged over 
all $N \times N$ matrices $U$ in the CUE, and when $k = 1, 2, 
3,$ and $4$ they obtained  the same values for the numbers corresponding 
to $g_k$ as those above. They argued that one could therefore 
model the moments $I_{k}(1/2, T)$ by the average of $|Z_{N}(U,\theta)|^{2k}$ 
over CUE and conjectured that 
$$
g_k = (k^{2}!)\prod_{j=0}^{k-1} \frac{j!}{(j+k)!}\,.
$$
Interestingly, the Keating--Snaith and Conrey--Gonek conjectures 
were first publicly announced at the Riemann Hypothesis Conference in 
Vienna, just moments after it was checked that the Keating--Snaith
conjecture in fact predicts that $g_{4} =24024$\,.

The characteristic polynomial model has proven to be extremely 
powerful for  predicting other behavior of the zeta--function  and L--functions 
that once seemed hopelessly beyond reach. In fact,
to a large extent it has been responsible for an explosion of 
activity in the field and 
of collaboration between number theorists and theoretical physicists.

Impressive as the characteristic polynomial model has proven to be, 
it has the obvious drawback that it contains no 
arithmetical information. The prime numbers do not appear in 
this model of the zeta--function! In the moment problem, this is 
reflected by the absence of the arithmetical factor $a_{k}$ in the Keating--Snaith 
conjecture. They had to insert it in an ad hoc way. Fortunately, 
in the moment problem, it was only the factor $g_{k}$ and not $a_{k}$ that 
proven elusive. A  precise and more satisfactory model for the zeta--function (and 
other L--functions) clearly has to include such relevant information.

In work in progress with J. Keating and  C. Hughes, we have now succeeded in finding 
such a model for $\zeta(s)$ and it can easily be generalized
to model any L--function.  I will conclude this lecture by describing 
the new model.

Roughly, we have proven that if the Riemann Hypothesis is true, then 
for $t\in \mathbb{R}$ and $ x > 1$ we have
$$
\zeta(\frac{1}{2} + i t) = \exp \left(\sum_{2 \leq  n \leq x}\,
\Lambda(n) / ({n}^{1/2 + i t} \log n) \right) 
\prod_{n}\exp \left(E_{1}((t-\gamma_n)\log x) \right) \\
$$
(times an error term that is essentially 1), where 
$E_{1}(z) = \int_{z}^{\infty} \frac{e^{-w}}{w} \,dw $ is the exponential integral.
I say ``roughly'' because one  also has to include smooth weights in the 
various factors. A similar formula holds throughout the critical 
strip. Since we expect the ordinates $\gamma_{n}$ of the zeros to behave like the 
eigenangles $\theta_{n}$ of $N \times N$ random matrices in CUE, 
and "scaling" suggests that we take $N$ to be the nearest integer to
$\frac{1}{2 \pi}\log t$, we take as our model for zeta 
$$
\exp \left(\sum_{2 \leq  n \leq x}\,
\Lambda(n) / ({n}^{1/2 + i t} \log n) \right) 
\prod_{n \leq N}\exp \left(E_{1}((\theta -\theta_n)\log x) \right) 
\,.
$$
The presence of the exponential integral makes it a little
complicated to compare this with the previous model, 
$$
Z_{N}(U,\theta) = \prod_{n=1}^{N} ( 1 - e^{i(\theta_{n} - \theta)} )\,.
$$
We note, however,that if $\theta_{n}$ is not too near $\theta$, 
then the new model looks approximately like
$$
\prod_{p \leq x} \left( 1 -  {p}^{-(1/2 + i t)} \right)^{-1} 
\prod_{n \leq N}\exp \left(1 - x^{ i (\theta -\theta_n)} \right) 
\,.
$$
Here we clearly see both the primes and the zeros, and how 
the parameter $x$ serves to connect them. The moments
$I_{k}(1/2, T)$ should now be given by the product of two moments--
one being the $2k$th power of the modulus of the product over primes 
integrated with respect to t,  the other being the $2k$th power of the modulus of
the product over the eigenangles averaged over the Circular Unitary Ensemble. 
We call the conjecture that the mean can be computed this way, that 
is, as a product of two different types of means, the ``Splitting Conjecture''.

The new model seems promising for many other investigations as well.
To give just one example, we hope to use it to understand the  horizontal 
distribution of the zeros of $\zeta^{\prime}(s)$ in the right half of the 
critical strip, a problem that has long defied us.
We also expect it to give us more insight into the connection between primes and 
zeros. If we are extremely lucky, perhaps we will even find explicit and useful connections  
between primes in special sequences, such as twin primes or primes of the 
form  $n^{2} + 1$, and the zeros. 

%%%%%%%%
%%%%%%%%%

\end{document}